\documentclass[12pt]{article}
\usepackage{amsxtra,amssymb,amsthm,amsmath,latexsym}

\theoremstyle{plain}
\def\R{{\mathbb R}}

\def\oH{{\overset{\circ}{H}}}
\def\oH1{{\overset{\circ}{H}\kern-.02in{}^1}}

\def\bee{\begin{equation*}}
\def\eee{\end{equation*}}
\def\be{\begin{equation}}
\def\ee{\end{equation}}

\begin{document} %\begin{titlepage}
\title{ Three open problems in analysis}

\author{A.G. Ramm\\
 Mathematics Department, Kansas State University, \\
 Manhattan, KS 66506-2602, USA\\
ramm@math.ksu.edu,\\ fax 785-532-0546, tel. 785-532-0580}
%http://www.math.ksu.edu/\,$\widetilde{\ }$\,ramm}

\date{}
\maketitle\thispagestyle{empty}

%\begin{abstract}
\footnote{MSC: 44A12, 35J55, 47B99     }
\footnote{key words: Radon transform, reaction-diffusion systems,
operator theory}

%\end{abstract}
%\end{titlepage}

\section{Statement of the problems }\label{S:1} 
In this
paper some open problems in analysis are formulated. These
problems were explained by the author at ICMAA6.

1. {\bf Injectivity of the classical Radon transform}.

Consider the Radon transform:
\be\label{e.1} Rf:=\int_{\ell_{\alpha, p}}fds, \ee
where  $\ell_{\alpha, p}$ is a straight line $\alpha \cdot
x=p$ on the plane $x=\{x_1, x_2\}$, $\alpha$ is a unit vector, $p$
is a real number, $ds$ is the element of the arclength of the straight
line.

Assume that 
\be\label{e.2} f\in L^1(\ell_{\alpha, p}) \ee 
for all
$p$   and $\alpha$, that $f$ is a continuous function, and that
\be\label{e.3} |f(x)|\leq c(1+|x|^m),\ee

where $c=const>0$, and $m\geq 0$ is a fixed number. Assume that
\be\label{e.4} Rf=0 \ee
for all $p$ and $\alpha$.

{\bf Problem 1:} {\it Does it follow from the above assumptions that 
$f=0$?}

There is a large literature on Radon transform
(see, e.g.,\,\cite{H}, \cite{R1} and references therein).
It is known (see, e.g,\, \cite{A}, \cite{R1})  that there are entire 
functions
not vanishing identically, such that (2) and (4) hold.

The open problem is to understand what the weakest natural restriction on
the growth of $f$ at infinity is for the Radon transform to be
injective. In other words, under what weakest growth restriction at
infinity do assumptions (1) and (3) imply $f=0$?

It is known (see \cite{R1}) that if $f\in L^1(\R^2, \frac 1{1+|x|})$
and (1) holds, then $f=0$, i.e., the Radon transform is injective on 
$ L^1(\R^2, \frac 1{1+|x|)})$.

2. {\bf A uniqueness problem}.

Let $L$ and $M$ be elliptic, second order, selfadjoint, strictly
positive Dirichlet operators in a bounded domain $D$ in $R^n, n>1$,
with a smooth connected boundary $S$, and the coefficients of $L$
and $M$ are real-valued functions, so that all the functions below
are real-valued. Let $a(x)$ and $b(x)$ be strictly positive functions, 
smooth in the closure of $D$. Let
\bee Lu+a(x)v=0\hbox{\ in\ } D,\quad -b(x)u+Mv=0\hbox{\ in\ }D,
\quad u=v=0\hbox{\ on\ }S.\eqno{(*)}\eee
{\bf Problem 2:} {\it Does $(*)$ imply}
\bee u=v=0\hbox{\  in\ } D?  \eqno{(**)} \eee
It is of no interest
to give sufficient conditions for $(**)$ to hold, such as, e.g.,  $|b-a|$ 
is
small, or $L=M$, or some other conditions. What is of interest is to
answer the question as stated, without any additional assumptions,
by either proving $(**)$ or constructing a counterexample.

In the one-dimensional case the answer to the question $(**)$ is yes
(see \cite{LG}).

3.{\bf A problem in operator theory}.

The question in two different forms is stated below as
Problem 3.1 and Problem 3.2.
These problems are closely related.

3.1.  Let $D$ be a bounded domain in $R^3$,
$D$ can be a box or a ball,
$f(x)\in L^2(D)$ be a function
 $f\not\equiv 0$, $F(z):=\int_Df(x) \exp(iz\cdot x)dx$, where
$z\in C^3$ is an
 entire function of exponential type.
 Let $L_j(z)$, $j=1,2$, be polynomials,
 $\text { deg} L_j(z)\geq 1$,
 $\mathcal L_j:=\{z: z\in C^3, L_j(z)=0 \}$
 be algebraic varieties.

 Define Hilbert spaces $H_j:=L^2(\mathcal L_j,dm_j)$, where
 $dm_j(z)$ are smooth rapidly decaying strictly positive measures
 on $\mathcal L_j$,
such that any exponential $exp(iz\cdot x)$ with any $x\in R^3$
belongs to $H_j$.
 Define linear operator $T$ from $H_1$ into $H_2$ by the formula:
 $Th:=\int_{\mathcal L_1}dm_1h(u_1)F(u_1+u_2):=g(u_2)$,
where $u_j\in \mathcal L_j$,
 $h\in H_1$, $g\in H_2$ for any
$h\in H_1$ (such are measures $dm_j$ by
the assumption).
 Assume that $\mathcal L_1$ and $\mathcal L_2$ are transversal,
which by
definition means that there exist
two points, one in $\mathcal L_1$
and one in  $\mathcal L_2$, such that the union of the bases of
 the tangent spaces to $\mathcal L_1$ and to $\mathcal L_2$
at these points, form a basis in $C^3$.
The same setting is of interest in dimension $n>3$ as well.

{\bf Problem 3.1}: {\it Is it true that T is not a finite-rank operator?}

{\it That is, is it true that dimension of the range of T is
infinite?}

{\bf Remark}: The assumption that $f(x)$ is in $L^2(D)$ is important:
if, for example, $f(x)$ is a delta function, then the answer to
the question of problem 3.1 is no:
the dimension of the range of $T$ in this case is one
if the delta function is supported at one point.

3.2. In the notations of Problem 3.1, choose points
$p_m\in \mathcal L_2$,
$m=1,2,\dots,M$, where $M$ is an arbitrary large fixed integer.
Consider the set $S$ of $M$ functions ${F(z+p_m)}$, $m=1,2,\dots,M$,
where $z\in \mathcal L_1$, and $F(z)$ is defined above:
it is
the Fourier transform of a compactly supported $L^2(D)$ function,
where $D$ is a bounded domain in $R^n$, $n>1$, $D$ can be 
a box or a ball.

{\bf Problem 3.2}: {\it Can one choose $p_m\in \mathcal L_2$ such that
the above set
$S$ of $M$ functions is linearly independent?}

In other words, can one choose $p_m\in \mathcal L_2$,
$m=1,2,\dots,M$,
such that the relation:
\bee
 \sum_{m=1}^M c_mF(z+p_m)=0 \quad\forall z\in \mathcal L_1
 \eqno{(***)}
\eee
implies $c_m=0$ for all $m=1,2,\dots,M$?
Here $c_m$ are constants.

These questions arise in the study of Property C (\cite{R2}).

\end{document}